\documentclass[11pt,a4paper]{amsart}%
\pdfoutput=1
\usepackage{hyperref}
\usepackage{amsfonts}
\usepackage{amsmath}
\usepackage{amssymb}
\usepackage{graphicx}
\usepackage{algorithm}
\usepackage{algorithmic}%
\newtheorem{theorem}{Theorem}
\newtheorem{lemma}[theorem]{Lemma}
\theoremstyle{definition}
\newtheorem{definition}[theorem]{Definition}
\newtheorem{proposition}[theorem]{Proposition}
\newtheorem{example}[theorem]{Example}

\newtheorem{corollary}[theorem]{Corollary}
\theoremstyle{remark}
\newtheorem{remark}[theorem]{Remark}
\numberwithin{equation}{section}
\theoremstyle{plain}

\begin{document}
\title[A framework for tropical mirror symmetry]{A framework for tropical mirror symmetry}
\author{Janko B\"{o}hm}
\address{Department of Mathematics, Universit\"{a}t des Saarlandes, Campus E2 4 \\
D-66123 \\
Saarbr\"{u}cken, Germany}
\email{boehm@math.uni-sb.de}
\thanks{The author was supported by DFG (German Research Foundation) through Grant BO3330/1-1.}

\begin{abstract}
Applying tropical geometry a framework for mirror symmetry including a mirror
construction for Calabi-Yau varieties was proposed by the author. We discuss
the conceptual foundations of this construction based on a natural mirror
map identifying deformations and divisors. We show how the construction
specializes to that by Batyrev for hypersurfaces and its generalization by
Batyrev and Borisov to complete intersections. Based on an explicit example we
comment on the implementation in the \textsc{Macaulay2} package
\textsc{SRdeformations}.

\end{abstract}
\maketitle

\section{Introduction\label{Sec Introduction}}

Mirror symmetry is a key link between mathematics and theoretical physics,
e.g., algebraic geometry obtains new ideas in enumerative geometry from
superstring theory which, in return, benefits from the study of Calabi-Yau
varieties. Important insight to mirror symmetry is gained by explicit
constructions computing for a given Calabi-Yau variety the corresponding
mirror Calabi-Yau, for a general account of the topic see \cite{CK Mirror
Symmetry and Algebraic Geometry}.

The mirror of the general quintic hypersurface in $\mathbb{P}^{4}$ was given
by Greene and Plesser \cite{GrPl Duality in CalabiYau moduli space} as an
orbifold of a $1$-parameter family of quintics. For toric hypersurfaces this
class of orbifolding constructions was unified by Batyrev \cite{Batyrev Dual
polyhedra and mirror symmetry for CalabiYau hypersurfaces in toric varieties}
using the involution of Gorenstein toric Fano varieties given by dualization
of reflexive polyhedra. Batyrev's description proved to be well suited for the
study of further properties like mirror duality of stringy Hodge numbers
\cite{Batyrev Stringy Hodge}, Picard-Fuchs equations \cite{Batyrev Variations
of the mixed Hodge structure of affine hypersurfaces in algebraic tori} and
much more. It was generalized by Batyrev and Borisov to complete intersections
\cite{BB Mirror duality and stringtheoretic Hodge numbers, BB On CalabiYau
complete intersections in toric varieties in HigherDimensional Complex
Varieties Trento 1994} using nef partitions and by Batyrev and Nill
\cite{BatyrevNill} via reflexive Gorenstein cones.

Based on ideas of Leung and Vafa \cite{LV Branes and Toric Geometry} and
Kontsevich and Soibelman \cite{KS Homological mirror symmetry and torus
fibrations}, Gross and Siebert \cite{GrSi Affine manifolds log structures and
mirror symmetry, GrSi Mirror symmetry via logarithmic degeneration data I,
GrSi From real affine geometry to complex geometry} used toric degenerations
and integrally affine manifolds to give a mirror construction, which is
expected to eventually relate $B$-model period integrals and tropically
counted $A$-model Gromov-Witten invariants. For a first instance of this in
the case of $\mathbb{P}^{2}$ see \cite{GrossP2}. For the fundamental idea of
tropical curve counting see Mikhalkin \cite{Mikhalkin}. In \cite{Gross} Gross
shows how to construct complete intersection mirrors. In order to apply the
Gross-Siebert program one has to obtain simple affine structures, which is
achieved by considering more general families, which become toric
degenerations after desingularization. For an independent construction of
these integral affine structures in the complete intersection case by Haase
see \cite{Haase}.

In \cite{B Mirror symmetry and tropical geomety} the author developed via
embedded tropical varieties a general framework for mirror symmetry leading to
an algorithmic mirror construction. It directly specializes to the known
constructions by Batyrev for hypersurfaces and its generalization by Batyrev
and Borisov to complete intersections, and reproduces that by R\o dland
\cite{Rodland The Pfaffian CalabiYau its Mirror} for a Pfaffian non-complete
intersection. The tropical mirror construction extends this construction to a
considerably larger class of Calabi-Yau varieties and produces explicit new
mirror examples \cite[Sec. 10.5]{B Mirror symmetry and tropical geomety}. It
comes with a natural mirror map identifying deformations and divisors. In this
paper we focus on the conceptual foundation of the tropical mirror
construction and how to recover the Batyrev-Borisov mirror of a complete
intersection. We also discuss the implementation in the \textsc{Macaulay2}
\cite{GS} package \textsc{SRdeformations} \cite{SRdeformations}.

\textbf{Acknowledgements. }I would like to thank Frank-Olaf Schreyer, Wolfram
Decker and Stavros Papadakis for many important suggestions while working on
my thesis and Victor Batyrev, Jan Arthur Christophersen, Mark Gross, Bernd
Siebert, Duco van Straten and Bernd Sturmfels for a lot of fruitful
discussions in the context of mirror symmetry.

\section{Ingredients from toric and tropical geometry, and deformation theory}

In this section we discuss basic facts from toric and tropical geometry and
deformation theory necessary to formulate the tropical mirror construction,
and fix some notation in this context.

\subsection{Toric geometry}

We introduce the basic toric objects used in the tropical mirror construction.
For more details on toric geometry see, e.g., \cite{CK Mirror Symmetry and
Algebraic Geometry} and \cite{Wisniewski Toric Mori theory and Fano manifolds}.

\subsubsection{Toric Fano varieties}

A \textbf{Fano polytope} $P\subset N_{\mathbb{R}}=N\otimes_{\mathbb{Z}%
}\mathbb{R}$ in the lattice $N=\mathbb{Z}^{n}$ is an integral polytope which
contains $0$ as its unique interior lattice point. The fan $\Sigma
=\Sigma\left(  P\right)  $ over $P$ given by the cones
$\operatorname*{posHull}F$ spanned by the faces $F$ of $P$ defines a
$\mathbb{Q}$-Gorenstein toric Fano variety $Y=\operatorname*{TV}\left(
\Sigma\right)  $ of dimension $n$. The toric strata of $Y$ correspond to the
faces of the dual polytope $\Delta=P^{\ast}\subset M_{\mathbb{R}}%
=M\otimes_{\mathbb{Z}}\mathbb{R}$ where $M=\operatorname*{Hom}\left(
N,\mathbb{Z}\right)  $.

Denote by $\mathbb{P}\left(  \Delta\right)  $ the projective toric variety
defined an integral polytope $\Delta\subset M_{\mathbb{R}}$ and by
$\Sigma=\operatorname{NF}\left(  \Delta\right)  =\Sigma\left(  P\right)  $ the
normal fan of $\Delta$, so $\mathbb{P}\left(  \Delta\right)  \cong%
\operatorname*{TV}\left(  \Sigma\right)  $ and $\mathcal{O}_{\mathbb{P}\left(
\Delta\right)  }\left(  1\right)  \cong\mathcal{O}_{\mathbb{P}\left(
\Delta\right)  }\left(  D_{\Delta}\right)  $ with the Cartier divisor%
\[
D_{\Delta}=%
{\textstyle\sum\nolimits_{r\in\Sigma\left(  1\right)  }}
-\min_{m\in\Delta}\left\langle m,\hat{r}\right\rangle D_{r}%
\]
Here we denote by $\Sigma\left(  1\right)  $ the set of rays (cones of
dimension $1$) of $\Sigma$, by $D_{r}$ the torus invariant prime Weil divisor
corresponding to $r\in\Sigma\left(  1\right)  $ and by $\hat{r}$ the minimal
lattice generator of $r$.

A polytope $\Delta\subset M_{\mathbb{R}}$ of dimension $n$ is called
\textbf{reflexive} if $\Delta$ and its dual $\Delta^{\ast}$ are integral and
contain $0$ in their interior. Then $\mathbb{P}\left(  \Delta\right)  $ and
$\mathbb{P}\left(  \Delta^{\ast}\right)  $ are Gorenstein toric Fano varieties.

\subsubsection{Cox ring}

Subvarieties of a toric variety $Y=\operatorname*{TV}\left(  \Sigma\right)  $
can be described by ideals in the \textbf{Cox ring} (or homogeneous coordinate
ring) $S=\mathbb{C}\left[  x_{r}\mid r\in\Sigma\left(  1\right)  \right]  $ of
$Y$, see \cite{Cox The homogeneous coordinate ring of a toric variety}. This
is a polynomial ring with one variable $x_{r}$ for each ray $r\in\Sigma\left(
1\right)  $ graded by the presentation sequence
\begin{equation}
0\rightarrow M\overset{A}{\rightarrow}\mathbb{Z}^{\Sigma\left(  1\right)
}\overset{\deg}{\rightarrow}A_{n-1}\left(  Y\right)  \rightarrow0
\label{Chow presentation}%
\end{equation}
of the Chow group $A_{n-1}\left(  Y\right)  $ of classes $\left[  D\right]  $
of Weil divisors $D$ modulo linear equivalence. The rows of the matrix $A$ are
given by the minimal lattice generators of the rays of $\Sigma$. In terms of
monomials%
\[
\deg\left(
{\textstyle\prod\nolimits_{r}}
x_{r}^{a_{r}}\right)  =\left[
{\textstyle\sum\nolimits_{r}}
a_{r}D_{r}\right]
\]

\begin{proposition}
The vector space of global sections of the reflexive sheaf of sections
$\mathcal{O}_{Y}\left(  D\right)  $ of a Weil divisor $D$ in $Y$ is isomorphic
to the degree $\left[  D\right]  $-part of the Cox ring%
\[
H^{0}\left(  Y,\mathcal{O}_{Y}\left(  D\right)  \right)  \cong S_{\left[
D\right]  }%
\]

\end{proposition}

We denote by $\Delta_{D}\subset M_{\mathbb{R}}$ the polytope of sections of a
divisor $D$, i.e., the convex hull of the torus invariant sections.

From $A_{n-1}\left(  Y\right)  $, depending only on the rays of $\Sigma$, and
the \textbf{irrelevant ideal}%
\[
B\left(  \Sigma\right)  =\left\langle
{\textstyle\prod\nolimits_{r\in\Sigma\left(  1\right)  ,\text{ }%
r\not \subset \sigma}}
x_{r}\mid\sigma\in\Sigma\right\rangle \subset S\text{,}%
\]
the toric variety $Y$ can be recovered as the categorical quotient%
\[
Y=\left(  \mathbb{C}^{\Sigma\left(  1\right)  }-V\left(  B\left(
\Sigma\right)  \right)  \right)  //\operatorname*{Hom}\nolimits_{\mathbb{Z}%
}\left(  A_{n-1}\left(  Y\right)  ,\mathbb{C}^{\ast}\right)
\]


\begin{definition}
If $I\subset S$ is generated by homogeneous elements $f\in S$ with
$\deg\left(  f\right)  \in\operatorname*{Pic}\left(  Y\right)  $, then $I$ is
\index{Pic-generated|textbf}%
called $\operatorname*{Pic}\left(  Y\right)  $\textbf{-generated}. The ideal
$I$ is called $\operatorname*{Pic}\left(  Y\right)  $\textbf{-saturated} if
$I_{\alpha}=\left(  I:B\left(  \Sigma\right)  ^{\infty}\right)  _{\alpha}$ for
all $\alpha\in\operatorname*{Pic}\left(  Y\right)  $.
\end{definition}

\begin{definition}
The \textbf{Picard-Cox ring} of $Y$ is%
\[
R=\bigoplus\nolimits_{\alpha\in\operatorname*{Pic}\left(  Y\right)  }%
S_{\alpha}%
\]

\end{definition}

By \cite{Cox The homogeneous coordinate ring of a toric variety}, if $Y$ is
simplicial, there is a one-to-one correspondence between the
$\operatorname*{Pic}\left(  Y\right)  $-generated and $\operatorname*{Pic}%
\left(  Y\right)  $-saturated ideals $I\subset S$ and the closed subschemes of
$Y$. Equivalently one can consider graded ideals of the Picard-Cox ring $R$,
which are saturated in $B\left(  \Sigma\right)  \cap R$.

\subsection{Tropical geometry\label{sec tropical varieties}}

Tropical geometry will be applied in the mirror construction as a tool to
explore one parameter degenerations with fibers in a toric variety, as it
associates to such a degeneration a combinatorial object. We recall some basic
facts, for more details on tropical geometry see, e.g., \cite{StML}.

\subsubsection{Amoebas}

Tropical geometry was motivated by the study of the \textbf{amoeba} of a
subvariety $V\subset\left(  \mathbb{C}^{\ast}\right)  ^{n}$ which is defined
as the image of $V$ under the map%
\begin{align*}
\log_{t}  &  :\left(  \mathbb{C}^{\ast}\right)  ^{n}\rightarrow\mathbb{R}%
^{n}\\
\left(  z_{1},...,z_{n}\right)   &  \mapsto\left(  \log_{t}\left\vert
z_{1}\right\vert ,...,\log_{t}\left\vert z_{n}\right\vert \right)
\end{align*}
for some base $t$. Note, that considering the fibers of this map relates
tropical geometry to the context of torus fibrations in mirror symmetry. The
limit of the amoeba for $t\rightarrow\infty$ in the Hausdorff metric on
compacts can be obtained as a non-Archimedian version of the amoeba:

\subsubsection{Tropical varieties\label{sub section tropical varieties}}

Consider the field of Puiseux series $\mathbb{C}\left\{  \left\{  t\right\}
\right\}  $, which is equipped with the valuation%
\begin{gather*}
val:\mathbb{C}\left\{  \left\{  t\right\}  \right\}  \rightarrow\mathbb{Q\cup
}\left\{  \infty\right\} \\
\sum_{j\in J}\alpha_{j}t^{j}\mapsto\min J
\end{gather*}
and with a norm $\left\Vert f\right\Vert =e^{-val\left(  f\right)  }$. Extend
$val$ and $\left\Vert -\right\Vert $ to the metric completion $K$ of
$\mathbb{C}\left\{  \left\{  t\right\}  \right\}  $ containing those elements
$\sum_{j\in J}\alpha_{j}t^{j}$, which satisfy the condition that any subset of
$J$ has a minimum. So $K$ is a complete algebraically closed non-Archimedian
field with surjective valuation $val:K\rightarrow\mathbb{R}\cup\left\{
\infty\right\}  $.

Let $I$ be an ideal in $K\left[  x_{1},...,x_{n}\right]  $. The image of the
algebraic variety $V_{K}\left(  I\right)  \subset\left(  K^{\ast}\right)
^{n}$ defined by $I$ under the non-Archimedian amoeba map%
\begin{gather*}
\operatorname*{val}=\log\left\Vert -\right\Vert :\left(  K^{\ast}\right)
^{n}\rightarrow\mathbb{R}^{n}\\
\left(  x_{1},...,x_{n}\right)  \mapsto\left(  val\left(  x_{1}\right)
,...,val\left(  x_{n}\right)  \right)
\end{gather*}
is called the non-Archimedian amoeba of $V_{K}\left(  I\right)  $ or
\textbf{tropical variety} $T\left(  I\right)  $ of $I$.

For $w\in\mathbb{R}^{n}$ the \textbf{initial form} $in_{w}\left(  f\right)  $
of $f\in K\left[  x_{1},...,x_{n}\right]  $ is the sum of the terms of maximal
weight with respect to $w$ and $weight\left(  c\right)  =-val\left(  c\right)
$ for $c\in K$. For any ideal $J\subset K\left[  x_{1},...,x_{n}\right]  $ its
\textbf{initial ideal} is%
\[
in_{w}\left(  J\right)  =\left\langle in_{w}\left(  f\right)  \mid f\in
J\right\rangle
\]

The \textbf{tropical semiring} is $\mathbb{R\cup}\left\{  \infty\right\}  $
with tropical addition and multiplication%
\begin{align*}
a\oplus b  &  =\min\left(  a,b\right) \\
a\odot b  &  =a+b
\end{align*}
For any polynomial%
\[
f=\sum\nolimits_{a}b_{a}\left(  t\right)  \cdot x^{a}\in K\left[
x_{1},...,x_{n}\right]
\]
define its \textbf{tropicalization} as the piecewise linear function%
\[
\operatorname*{trop}\left(  f\right)  =%
{\textstyle\bigoplus\nolimits_{a}}
val\left(  b_{a}\left(  t\right)  \right)  \odot x^{\odot a}%
\]
and by $T\left(  \operatorname*{trop}\left(  f\right)  \right)  $ its corner
locus, i.e., the set of $w\in\mathbb{R}^{n}$ such that the minimum is attained
at least twice. Then the fundamental theorem of tropical geometry is:

\begin{theorem}
\cite[Sec. 9.2]{Sturmfels Solving Systems of Polynomial Equations}, \cite[Sec.
2]{SpS The tropical Grassmannian},\cite{StML} If $I\subset K\left[
x_{1},...,x_{n}\right]  $ is an ideal, then%
\begin{align*}
T\left(  I\right)   &  =\left\{  w\in\mathbb{R}^{n}\mid in_{w}\left(
I\right)  \text{ contains no monomial}\right\} \\
&  =%
{\textstyle\bigcap\nolimits_{f\in I}}
T\left(  \operatorname*{trop}\left(  f\right)  \right)
\end{align*}

\end{theorem}

\begin{remark}
The tropical variety $T\left(  I\right)  $ has a structure of a polyhedral
cell complex, its dimension is the Krull dimension of $K\left[  x_{1}%
,...,x_{n}\right]  /I$ and it is equidimensional if $V_{K}\left(  I\right)  $ is.
\end{remark}

The tropical variety $T\left(  I\right)  $ is a subset of the space of weight
vectors $\left(  w_{x_{1}},...,w_{x_{n}}\right)  \in\mathbb{R}^{n}$ on the
monomials of $K\left[  x_{1},...,x_{n}\right]  $.

\subsection{Bergman fan\label{sec Bergman fan}}

From the point of view of Gr\"{o}bner fans, see for example \cite{Sturmfels
Groebner Bases and Convex Polytopes}, it is more natural to work with a
tropical fan: We give a non-Archimedian definition of the Bergman fan
(historically defined via an amoeba type limit). Let $I\subset\mathbb{C}%
\left[  t\right]  \left[  x_{1},...,x_{n}\right]  $ be an ideal and denote by
$L$ the metric completion of $\mathbb{C}\left\{  \left\{  s\right\}  \right\}
$. The image of $V_{L}\left(  I\right)  $ under%
\begin{gather*}
\left(  L^{\ast}\right)  ^{n+1}\rightarrow\mathbb{R}^{n+1}\\
\left(  t,x_{1},...,x_{n}\right)  \mapsto\left(  val\left(  t\right)
,val\left(  x_{1}\right)  ,...,val\left(  x_{n}\right)  \right)
\end{gather*}
is a fan (the tropical variety of $I$ considering $t$ as a variable), which we
denote as the \textbf{Bergman fan }of $I$.

For a fan $\Sigma$ and a hyperplane $H\ $in $\mathbb{R}^{n+1}$ define
$\Sigma\cap H$ as the polyhedral cell complex consisting of the faces
$\sigma\cap H$ for $\sigma\in\Sigma$. With this notation we immediately get:

\begin{proposition}
For an ideal $I\subset\mathbb{C}\left[  t\right]  \left[  x_{1},...,x_{n}%
\right]  $%
\[
\operatorname*{BF}\left(  I\right)  \cap\left\{  w_{t}=1\right\}  =T\left(
I\right)
\]

\end{proposition}

The intersection with the hyperplane $\left\{  w_{t}=1\right\}  $ amounts to
identification of the parameter $s$ of the Puiseux series solutions and the
parameter $t$ of the degeneration.

\subsection{Tropical varieties and the Cox ring}

Consider an ideal $I\subset K\otimes S$ where $S$ is the Cox ring of a toric
variety $Y$. As seen in Section \ref{sub section tropical varieties}, the
tropical variety of $I$ should be considered as a subset of the space of
weight vectors. Hence in the Cox setup $T\left(  I\right)  $ is naturally a
subset of the weight space%
\[
\frac{\operatorname*{Hom}\nolimits_{\mathbb{R}}\left(  \mathbb{R}%
^{\Sigma\left(  1\right)  },\mathbb{R}\right)  }{\operatorname*{Hom}%
\nolimits_{\mathbb{R}}\left(  A_{n-1}\left(  Y\right)  \otimes\mathbb{R}%
,\mathbb{R}\right)  }\underset{\varphi}{\overset{\cong}{\rightleftarrows}%
}N_{\mathbb{R}}%
\]
on $S$, where the isomorphism $\varphi$ is obtained from the presentation
sequence of the Chow group (Equation \ref{Chow presentation}).

For example if $Y=\mathbb{P}^{n}$ the space of weights will be $\mathbb{R}%
^{n+1}/\mathbb{R}\left(  1,...,1\right)  $.

\begin{remark}
When discarting the grading, the tropical variety may still contain linear
space after dividing by $\operatorname*{Hom}\nolimits_{\mathbb{R}}\left(
A_{n-1}\left(  Y\right)  \otimes\mathbb{R},\mathbb{R}\right)  $. In some
settings it makes sense to divide out also this lineality space, e.g., in the
context of tropical Grassmannians \cite{SpS The tropical Grassmannian}. In
general however, it should be considered as part of the tropical variety, as
then the dimensions of the tropical variety and the algebraic subvariety of
$Y$ will coincide.
\end{remark}

\subsection{Deformations of monomial ideals\label{sec def}}

Let $I_{0}$ be a reduced monomial ideal in the Cox ring $S$ of the toric
variety $Y$. As $I_{0}$ is generated by finitely many elements and the space
of elements of $S$ of a given degree is finite-dimensional, the degree $0$
homomorphisms in $\operatorname*{Hom}\left(  I_{0},S/I_{0}\right)  $ form a
finite-dimensional vector space denoted by $\operatorname*{Hom}\left(
I_{0},S/I_{0}\right)  _{0}$. The big torus $\left(  \mathbb{C}^{\ast}\right)
^{\Sigma\left(  1\right)  }$ acts by%
\[%
\begin{tabular}
[c]{lll}%
$\operatorname*{Hom}\nolimits_{\mathbb{Z}}\left(  \mathbb{Z}^{\Sigma\left(
1\right)  },\mathbb{C}^{\ast}\right)  \times\mathbb{C}\left[  \mathbb{Z}%
^{\Sigma\left(  1\right)  }\right]  $ & $\rightarrow$ & $\mathbb{C}\left[
\mathbb{Z}^{\Sigma\left(  1\right)  }\right]  $\\
\multicolumn{1}{c}{$\left(  \lambda,m\right)  $} & $\mapsto$ & $\lambda\left(
m\right)  \cdot m$%
\end{tabular}
\ \ \
\]
on $\mathbb{C}\left[  \mathbb{Z}^{\Sigma\left(  1\right)  }\right]  $ and on
$S$. The induced action of the abelian group $\operatorname*{Hom}%
\nolimits_{\mathbb{Z}}\left(  \mathbb{Z}^{\Sigma\left(  1\right)  }%
,\mathbb{C}^{\ast}\right)  $ on the vector space $\operatorname*{Hom}\left(
I_{0},S/I_{0}\right)  _{0}$ gives a representation
\[
\operatorname*{Hom}\nolimits_{\mathbb{Z}}(\mathbb{Z}^{\Sigma\left(  1\right)
},\mathbb{C}^{\ast})\rightarrow\operatorname*{GL}\left(  \operatorname*{Hom}%
\left(  I_{0},S/I_{0}\right)  _{0}\right)
\]
which decomposes into characters, as any irreducible representation of an
abelian group over an algebraically closed field is $1$-dimensional. So the
vector space $\operatorname*{Hom}\left(  I_{0},S/I_{0}\right)  _{0}$ has a
basis of deformations which are characters. Being of degree $0$, any such
homomorphism corresponds to an element $\alpha\in M\cong\operatorname*{image}%
\left(  A\right)  \subset\mathbb{Z}^{\Sigma\left(  1\right)  }$ (using the
notation of the presentation of the Chow group of divisors of $Y$ from
Equation \ref{Chow presentation}). On the other hand to any $\alpha\in M$ we
can associate a homomorphism $\delta_{\alpha}:I_{0}\rightarrow S/I_{0}$ (which
may be $0$). See also the results in \cite{AC Cotangent cohomology of
StanleyReisner rings}.

\section{Basic formulation of the tropical mirror
construction\label{sec tropical mirror}}

\subsection{Conceptual foundation}

Mirror symmetry is usually considered in the context of Calabi-Yau varieties,
i.e., normal projective algebraic variety with at worst Gorenstein canonical
singularities, trivial canonical sheaf $K_{X}=\Omega_{X}^{d}\cong%
\mathcal{O}_{X}$ and $h^{i}\left(  X,\mathcal{O}_{X}\right)  =0$ for $0<i<d$.
We will also consider this setup, however the tropical mirror construction is
not limited to the Calabi-Yau case.

\subsubsection{Complex and K\"{a}hler moduli}

In physics (see also \cite{CK Mirror Symmetry and Algebraic Geometry}) a
superconformal field theory is associated to a tuple $\left(  X,\omega\right)
$ of a Calabi-Yau variety $X$ and a complexified K\"{a}hler form $\omega=B+iJ$
with $B,J\in H^{2}\left(  X,\mathbb{R}\right)  $ and $J$ a K\"{a}hler class.
Mirror symmetry postulates the existence of a mirror dual tuple $\left(
X^{\circ},\omega^{\circ}\right)  $ leading to an isomorphic superconformal
field theory.

Keeping $\omega\ $fixed and varying $X$ should translate into $X^{\circ}$
being fixed and $\omega^{\circ}$ vary, and vice versa, the identification
given by the so called mirror map. Hence locally the complex moduli space of
$X$ is being identified with the K\"{a}hler moduli space of $\omega^{\circ}$.
So the corresponding tangent spaces $H^{1}\left(  T_{X}\right)  =H^{d-1,1}%
\left(  X\right)  $ of the complex moduli space and $H^{1,1}\left(  X^{\circ
}\right)  $ of the K\"{a}hler moduli space are isomorphic.

\subsubsection{Mirror symmetry and degenerations}

By this argument, mirror symmetry should be considered in a natural way not as
a relation on individual Calabi-Yau varieties, but rather on embedded flat
families. This idea is already present, e.g., in the representation of the
mirror of the general quintic $X$ as a $h^{2,1}\left(  X^{\circ}\right)
=h^{1,1}\left(  X\right)  =1$-parameter family degenerating in the union of
$5$ planes \cite{GrPl Duality in CalabiYau moduli space}. It was formalized in
\cite{GrSi Affine manifolds log structures and mirror symmetry, GrSi Mirror
symmetry via logarithmic degeneration data I, GrSi From real affine geometry
to complex geometry} in the context of toric degenerations. However note, that
some degenerations, one would like to apply mirror symmetry to (e.g., some
Pfaffian examples \cite{Rodland The Pfaffian CalabiYau its Mirror}, \cite{B
Mirror symmetry and tropical geomety}) do not fall in this category.

\subsubsection{Basic setup}

In the approach presented here we will see mirror symmetry as a correspondence
of monomial degenerations of Calabi-Yau varieties. So we consider flat
families $\mathfrak{X}\subset Y\times\operatorname*{Spec}\mathbb{C}\left[
t\right]  $ with Calabi-Yau fibers $X_{t}\subset Y$ in a $\mathbb{Q}%
$-Gorenstein toric Fano variety $Y$. The degeneration $\mathfrak{X}$ is
specified by an ideal $I\subset\mathbb{C}\left[  t\right]  \otimes S$,
homogeneous with respect to the variables of the Cox ring $S$ of $Y$, and
$X_{0}$ by a monomial ideal $I_{0}\subset S$.

The goal is to associate to $\mathfrak{X}$ a mirror degeneration
$\mathfrak{X}^{\circ}$ with fibers in a mirror toric Fano variety $Y^{\circ}$.
This will be done in a way, that we obtain a natural mirror map relating
$H^{d-1,1}\left(  X\right)  $ and $H^{1,1}\left(  X^{\circ}\right)  $ for
generic fibers $X\ $of $\mathfrak{X}$ and $X^{\circ}$ of $\mathfrak{X}^{\circ
}$. This can be seen as a generalization of the monomial-divisor mirror map
introduced in \cite{AGM The MonomialDivisor Mirror Map} for hypersurfaces in
Gorenstein toric Fano varieties.

For the construction we represent the complex moduli space of $X$ via a one
paramter family $\mathfrak{X}$ which is general in the following sense:
Consider a big torus invariant basis $v_{1},...,v_{p}\in\operatorname*{Hom}%
\left(  I_{0},S/I_{0}\right)  _{0}$ of degree $0$ homomorphisms of the tangent
space of the component of the Hilbert scheme of $X_{0}$ containing
$\mathfrak{X}$. Suppose that the tangent vector $v=\sum_{i=1}^{p}\lambda
_{i}v_{i}$ of $\mathfrak{X}$ satisfies $\lambda_{i}\neq0$ $\forall i$.

\begin{remark}
One could also consider special subfamilies, e.g., with prescribed
singularities. Furthermore, it seems possible to formulate a version of the
construction in several parameters $t_{1},...,t_{p}$, which could also be
handled by tropical geometry like in Section \ref{sec Bergman fan}. This may
avoid representing the moduli space by a one parameter family and eventually
could increase the scope of the construction.
\end{remark}

\subsubsection{Basic idea}

As discussed in Section \ref{sec def} the elements $v_{1},...,v_{p}$
correspond to elements $\alpha_{1},...,\alpha_{p}\in M$ of the lattice of
monomials of $Y$. The basic idea of the tropical mirror construction is to
consider the convex hull $\nabla^{\ast}$ of $\alpha_{1},...,\alpha_{p}$ and as
$Y^{\circ}$ the toric variety defined by the fan $\Sigma^{\circ}$ over the
faces of $\nabla^{\ast}$. Hence toric divisors of $Y^{\circ}$, and the induced
divisors on a prospective mirror inside constructed via tropical geometry,
will correspond to deformations of $X_{0}$. On the other hand the deformations
of the mirror special fiber $X_{0}^{\circ}$ should be induced by the toric
divisors of $Y$.

We now give a short outline of the general tropical mirror construction, for
more details see \cite{B Mirror symmetry and tropical geomety}.

\subsection{Input data}

We begin by summarizing the input data:

\subsubsection{Toric Fano variety}

Let $N=\mathbb{Z}^{n}$ and $\Delta^{\ast}\subset N_{\mathbb{R}}$ be a Fano
polytope and $Y=X\left(  \Sigma\right)  $, $\Sigma=\operatorname{Fan}\left(
\Delta^{\ast}\right)  $, the corresponding toric Fano variety with Cox ring
$S=\mathbb{C}\left[  x_{r}\mid r\in\Sigma\left(  1\right)  \right]  $ graded
by%
\[
0\rightarrow M\overset{A}{\rightarrow}\mathbb{Z}^{\Sigma\left(  1\right)
}\overset{\deg}{\rightarrow}A_{n-1}\left(  Y\right)  \rightarrow0
\]

\subsubsection{Monomial Calabi-Yau}

Let $X_{0}\subset Y$ be given by a reduced, $\operatorname{Pic}\left(
Y\right)  $-generated monomial ideal $I_{0}\subset S$ such that the subcomplex
$\operatorname*{Strata}\nolimits_{\Delta}\left(  I_{0}\right)  \subset
\partial\Delta$ of the boundary complex of $\Delta$, consisting of the toric
strata of $X_{0}$, is homeomorphic to a sphere.

\subsubsection{Degeneration}

Let $\mathfrak{X}\subset Y\times\operatorname*{Spec}\mathbb{C}\left[
t\right]  $ be a flat family of Calabi-Yau varieties of dimension $d$ with
fibers $X_{t}\subset Y$ and monomial special fiber $X_{0}$. The degeneration
$\mathfrak{X}$ is specified by an ideal $I\subset\mathbb{C}\left[  t\right]
\otimes S$, which is homogeneous with respect to the variables of $S$.


\subsection{Construction of the mirror polarization via Gr\"{o}bner
bases\label{sec tropical mirror polarization}}

Fix a monomial ordering $>$ on $\mathbb{C}\left[  t\right]  \otimes S$, which
is respecting the Cox grading on $S$ and is local in $t$, and denote by
$>_{w}$ the weight ordering by $w$ refined by $>$.

\begin{definition}
The \textbf{Gr\"{o}bner cone of special fiber weights} is defined as the
closed cone%
\[
C_{I_{0}}\left(  I\right)  =\left\{  -\left(  w_{t},w_{x}\right)
\in\mathbb{R}\oplus N_{\mathbb{R}}\mid L_{>_{(w_{t},\varphi(w_{x}))}}\left(
I\right)  =I_{0}\right\}
\]
with $\varphi$ as given in Section \ref{sec tropical varieties}.
\end{definition}

Under suitable conditions on the genericity of the degeneration and on the
smoothness of the base, the intersection of $C_{I_{0}}\left(  I\right)  $ with
the hyperplane of $t$-weight one is a polytope and its dual a Fano polytope.
We restrict to this case.

\begin{definition}
The \textbf{polytope of special fiber weights} is%
\[
\nabla_{I_{0}}\left(  I\right)  =C_{I_{0}}\left(  I\right)  \cap\left\{
w_{t}=1\right\}  \subset N_{\mathbb{R}}%
\]

\end{definition}

The Fano polytope $\nabla_{I_{0}}\left(  I\right)  ^{\ast}\subset
M_{\mathbb{R}}$ defines a toric Fano variety $Y^{\circ}=X\left(  \Sigma
^{\circ}\right)  $ of the same dimension as $Y$ by the fan $\Sigma^{\circ
}=\operatorname*{Fan}\left(  \nabla_{I_{0}}\left(  I\right)  ^{\ast}\right)
\subset M_{\mathbb{R}}$ over the faces of $\nabla_{I_{0}}\left(  I\right)
^{\ast}$. Denote by $S^{\circ}=\mathbb{C}\left[  z_{r}\mid r\in\Sigma^{\circ
}\left(  1\right)  \right]  $ the Cox ring of $Y^{\circ}$, graded by%
\[
0\rightarrow N\overset{A^{\circ}}{\rightarrow}\mathbb{Z}^{\Sigma^{\circ
}\left(  1\right)  }\overset{\deg}{\rightarrow}A_{n-1}\left(  Y^{\circ
}\right)  \rightarrow0
\]

\subsection{Tropical geometry construction of the mirror degeneration}

\subsubsection{Mirror special fiber}

Denote by $\partial C_{I_{0}}\left(  I\right)  $ the fan of all boundary faces
of the cone $C_{I_{0}}\left(  I\right)  $.

\begin{definition}
Consider the fan $BF_{I_{0}}\left(  I\right)  =BF\left(  I\right)
\cap\partial C_{I_{0}}\left(  I\right)  $ of the tropical faces of $C_{I_{0}%
}\left(  I\right)  $. By intersecting all cones of $BF_{I_{0}}\left(
I\right)  $ with the hyperplane $\left\{  w_{t}=1\right\}  $ one obtains a
subcomplex $T_{I_{0}}\left(  I\right)  \subset\partial\nabla_{I_{0}}\left(
I\right)  $, which we will denote as the \textbf{special fiber tropical
variety}.
\end{definition}

The support of $T_{I_{0}}\left(  I\right)  $ is a subset of the tropical
variety of $I$. The complex $T_{I_{0}}\left(  I\right)  $ is a subdivision of
the dual sphere of $\operatorname*{Strata}\nolimits_{\Delta}\left(
I_{0}\right)  $.

Recall, that there is a one-to-one correspondence between the Cox variables
$y_{r}$ of $Y^{\circ}$, the rays $r\in\Sigma^{\circ}\left(  1\right)  $ and
the facets $F_{r}^{\circ}$ of $\nabla$. Associated to the complex $T_{I_{0}%
}\left(  I\right)  \subset\partial\nabla$ of dimension $d$ we have a reduced
monomial ideal%
\begin{align}
I_{0}^{\circ}  &  =\left\langle
{\textstyle\prod\limits_{r\in J}}
y_{r}\mid J\subset\Sigma^{\circ}\left(  1\right)  \text{, }%
\operatorname*{supp}\left(  T_{I_{0}}\left(  I\right)  \right)  \subset%
{\textstyle\bigcup\limits_{r\in J}}
F_{r}^{\circ}\right\rangle \label{equ I0mirror}\\
&  =%
{\textstyle\bigcap\nolimits_{F^{\circ}\text{ facet of }T_{I_{0}}\left(
I\right)  }}
\left\langle y_{G^{\ast}}\mid G\text{ a facet of }\nabla\text{ with }F^{\circ
}\subset G\right\rangle \nonumber
\end{align}
defining a monomial Calabi-Yau of equi-dimension $d$ in $Y^{\circ}$. The ideal
is $\Sigma$-saturated, i.e., all primary components are strata of $Y^{\circ}$.
The first line says that $I_{0}^{\circ}$ is generated by the products of
variables which, seen as a union of facets of $\nabla_{I_{0}}\left(  I\right)
$, geometrically contain the support of $T_{I_{0}}\left(  I\right)  $. The
second line gives the unique irreducible decomposition of $I_{0}^{\circ}$. The
ideal of a maximal stratum $F^{\circ}\in T_{I_{0}}\left(  I\right)  $ is
generated by all variables which, considered as facets of $\nabla_{I_{0}%
}\left(  I\right)  $, contain $F^{\circ}$.

\subsubsection{First order mirror
degeneration\label{sec first order mirror degeneration}}

For a subcomplex $C$ of complex of faces of $\Delta$ denote by $C^{\ast
}\subset\Delta^{\ast}$ the co-complex of dual faces $F^{\ast}$ with $F\in C$.
We consider $\operatorname*{Strata}\nolimits_{\Delta}\left(  I_{0}\right)
^{\ast}\subset\Delta^{\ast}$ as the co-complex of deformations of the mirror.
The lattice points $\alpha\in\Xi=\operatorname*{supp}\left(
\operatorname*{Strata}\nolimits_{\Delta}\left(  I_{0}\right)  ^{\ast}\right)
\cap N$ of the support of this complex correspond via%
\[
0\rightarrow N\overset{A^{\circ}}{\rightarrow}\mathbb{Z}^{\Sigma^{\circ
}\left(  1\right)  }\overset{\deg}{\rightarrow}A_{n-1}\left(  Y^{\circ
}\right)  \rightarrow0
\]
to degree $0$ Cox Laurent monomials and represent degree $0$ deformations
$\varphi_{\alpha}\in\operatorname*{Hom}\nolimits_{S^{\circ}}\left(
I_{0}^{\circ},S^{\circ}/I_{0}^{\circ}\right)  _{0}$. Denote by $R^{\circ
}\subset S^{\circ}$ the Picard-Cox ring of $Y^{\circ}$.

\begin{definition}
\label{def tropmirror}The \textbf{first order tropical mirror} $\mathfrak{X}%
^{\circ}\subset Y^{\circ}\times\operatorname*{Spec}\mathbb{C}\left[  s\right]
/\langle s^{2}\rangle$ of $\mathfrak{X}$ is defined by the ideal%
\[%
\begin{tabular}
[c]{l}%
$I^{\circ}=\langle m^{\circ}+s\cdot%
{\textstyle\sum\nolimits_{\alpha\in\Xi}}
c_{\alpha}\cdot\varphi_{\alpha}(m^{\circ})\mid m^{\circ}\in I_{0}^{\circ}\cap
R^{\circ}\rangle\subset\mathbb{C}[s]/\langle s^{2}\rangle\otimes S^{\circ}$%
\end{tabular}
\ \
\]
with generic coefficients $c_{\alpha}$.
\end{definition}

Note, that it is sufficient to know a given family up to first order in the
case of complete intersections (due to the Koszul complex resolution) and
codimension $3$ Gorenstein varieties (due to the theorem of Buchsbaum and
Eisenbud, \cite{BE}).

\section{Application to Gorenstein complete
intersections\label{Sec Applications}}

\subsection{Setup\label{Sec Setup}}

We consider the setup of the mirror construction by Batyrev and Borisov
\cite{BB On CalabiYau complete intersections in toric varieties in
HigherDimensional Complex Varieties Trento 1994} for complete intersections in
Gorenstein toric Fano varieties. Let $Y=\mathbb{P}\left(  \Delta\right)  $ be
a Gorenstein toric Fano variety of dimension $n$, represented by the reflexive
polytope $\Delta\subset M_{\mathbb{R}}$, with normal fan $\Sigma\subset
N_{\mathbb{R}}$ and Cox ring $S$. A disjoint union%
\[
\Sigma\left(  1\right)  =J_{1}\cup...\cup J_{c}%
\]
is called a \textbf{nef partition} if all $E_{j}=\sum_{r\in J_{j}}D_{v}$ are
Cartier, spanned by global sections. By $\sum_{j=1}^{c}E_{j}=\sum_{r\in
\Sigma\left(  1\right)  }D_{r}=-K_{Y}$ general sections of $\mathcal{O}\left(
E_{1}\right)  ,...,\mathcal{O}\left(  E_{c}\right)  $ give a Calabi-Yau
complete intersection $X\subset Y$.

\subsection{Outline of the construction by Batyrev and
Borisov\label{sec BatyrevBorisov}}

For the setup from Section \ref{Sec Setup} Batyrev and Borisov construct the
mirror of $X$.

\begin{proposition}
\cite{BB On CalabiYau complete intersections in toric varieties in
HigherDimensional Complex Varieties Trento 1994} The polytopes $\Delta
_{j}=\Delta_{E_{j}}$ of sections of $E_{j}$ are lattice polytopes, and it
holds
\[
\Delta=\Delta_{1}+...+\Delta_{c}%
\]

\end{proposition}

Define the lattice polytope $\nabla_{j}$ as the convex hull%
\[
\nabla_{j}=\operatorname{convHull}\left\{  \left\{  0\right\}  \cup
J_{j}\right\}
\]
and $\nabla$ by%
\[
\nabla^{\ast}=\operatorname{convHull}\left(  \Delta_{1}\cup...\cup\Delta
_{c}\right)
\]

\begin{proposition}
\cite{BB On CalabiYau complete intersections in toric varieties in
HigherDimensional Complex Varieties Trento 1994} It holds $\nabla=\nabla
_{1}+...+\nabla_{c}$.
\end{proposition}

In particular $\nabla$ is a lattice polytope containing $0$, hence:

\begin{corollary}
\cite{BB On CalabiYau complete intersections in toric varieties in
HigherDimensional Complex Varieties Trento 1994} The polytope $\nabla$ is reflexive.
\end{corollary}

Let $\mathbb{P}\left(  \nabla\right)  $ be the Gorenstein toric Fano variety
associated to $\nabla$. Then%
\[
\sum_{j=1}^{c}D_{\nabla_{j}}=-K_{\mathbb{P}\left(  \nabla\right)  }%
\]
is a nef partition, and $X^{\circ}$ given by general sections of
$\mathcal{O}\left(  D_{\nabla_{1}}\right)  ,...,\mathcal{O}\left(
D_{\nabla_{c}}\right)  $ is a Calabi-Yau complete intersection in
$\mathbb{P}\left(  \nabla\right)  $.

\begin{theorem}
\cite{BB Mirror duality and stringtheoretic Hodge numbers} The Calabi-Yau
complete intersections $X$ and $X^{\circ}$ form a stringy topological mirror pair.
\end{theorem}

A maximal projective subdivision $\bar{\Sigma}$ of $\Sigma=\operatorname*{NF}%
\left(  \Delta\right)  $ gives a maximal projective partial crepant
desingularization%
\[
f:X\left(  \bar{\Sigma}\right)  \rightarrow\mathbb{P}\left(  \Delta\right)
\]
such that the $T$-divisors of the projective toric variety $X\left(
\bar{\Sigma}\right)  $ correspond to the lattice points of the boundary of
$\Delta^{\ast}$. Then $f$ induces a resolution $\bar{X}\rightarrow X$ of the
complete intersection $X\subset\mathbb{P}\left(  \Delta\right)  $ such that
$\bar{X}$ is a complete intersection, has at most Gorenstein terminal abelian
quotient singularities and $K_{\bar{X}}=\mathcal{O}_{\bar{X}}$. In particular,
if $\dim\left(  \bar{X}\right)  \leq3$, then $\bar{X}$ is smooth.

\subsection{Degenerations associated to complete
intersections\label{sec degen ci}}

The general nef complete intersection has a natural monomial degeneration
using the Koszul complex resolution:

\begin{lemma}
\label{lem degen ci}\cite{B Mirror symmetry and tropical geomety} Consider a
nef partition $\Sigma\left(  1\right)  =J_{1}\cup...\cup J_{c}$ as in the
setup of Section \ref{Sec Setup},%
\[
m_{j}=%
{\textstyle\prod\nolimits_{v\in J_{j}}}
x_{v}\in S
\]
and the reduced $\operatorname*{Pic}\left(  Y\right)  $-generated monomial
ideal
\[
I_{0}=\left\langle m_{j}\mid j=1,...,c\right\rangle
\]
Let $g_{j}\in S_{\left[  E_{j}\right]  }$ be general sections of
$\mathcal{O}\left(  E_{j}\right)  $ (corresponding to a general linear
combination of the lattice points of $\Delta_{E_{j}}$) not involving monomials
in $I_{0}$. Then the $\operatorname*{Pic}\left(  Y\right)  $-generated ideal%
\begin{equation}
I=\left\langle f_{j}=t\cdot g_{j}+m_{j}\mid j=1,...,c\right\rangle
\subset\mathbb{C}\left[  t\right]  \otimes S\nonumber
\end{equation}
defines a flat family $\mathfrak{X}\subset Y\times\operatorname*{Spec}\left(
\mathbb{C}\left[  t\right]  \right)  $ with fibers in $Y$ and special fiber
given by $I_{0}$.
\end{lemma}

The deformations of $I_{0}$ are unobstructed and the base space is smooth. Let
$v_{1},...,v_{p}\in\operatorname*{Hom}\left(  I_{0},S/I_{0}\right)  _{0}$ be a
basis of the tangent space of the Hilbert scheme of $X_{0}$. The degeneration
$\mathfrak{X}$ is general in the sense that if $v$ is the tangent vector of
$\mathfrak{X}$ and $v=\sum_{i=1}^{p}\lambda_{i}v_{i}$, then we have
$\lambda_{i}\neq0$ $\forall i$.

\subsection{Tropical construction of the Batyrev-Borisov mirror}

We now apply the tropical mirror construction to the canonical degeneration
$\mathfrak{X}$ from Section \ref{sec degen ci} of a given nef complete
intersection, and show that $\mathfrak{X}^{\circ}$ is the canonical
degeneration associated to the Batyrev-Borisov mirror. For details see \cite{B
Mirror symmetry and tropical geomety}.

\subsubsection{Construction of the mirror polarization via Gr\"{o}bner bases
techniques}

As $f_{1},...,f_{c}$ form a reduced Gr\"{o}bner basis with respect to any
monomial ordering selecting $I_{0}$ as lead ideal, we get:

\begin{lemma}
For the degeneration defined in Lemma \ref{lem degen ci} the special fiber
Gr\"{o}bner cone is%
\[
C_{I_{0}}\left(  I\right)  =\left\{  \left(  w_{t},w_{x}\right)  \in
\mathbb{R}\oplus N_{\mathbb{R}}\mid\operatorname*{trop}\left(  g_{j}\right)
\left(  \varphi\left(  w_{x}\right)  \right)  +w_{t}\geq\operatorname*{trop}%
\left(  m_{j}\right)  \left(  \varphi\left(  w_{x}\right)  \right)  \text{
}\forall j\right\}
\]

\end{lemma}

\noindent Hence the dual of the special fiber polytope
\[
\nabla_{I_{0}}\left(  I\right)  =\left\{  w_{x}\in N_{\mathbb{R}}%
\mid\operatorname*{trop}\left(  g_{j}\right)  \left(  \varphi\left(
w_{x}\right)  \right)  +1\geq\operatorname*{trop}\left(  m_{j}\right)  \left(
\varphi\left(  w_{x}\right)  \right)  \text{ }\forall j\right\}
\]
can be described as the convex hull of the lattice monomials $A^{-1}(\frac
{m}{m_{j}})$ with the monomials $m$ appearing in $g_{j}$ for $j=1,...,c$, hence:

\begin{corollary}
With the notation from Section \ref{sec BatyrevBorisov} and the degeneration
defined in Lemma \ref{lem degen ci}
\[
\nabla_{I_{0}}\left(  I\right)  ^{\ast}=\operatorname{convHull}\left(
\Delta_{1}\cup...\cup\Delta_{c}\right)  =\nabla^{\ast}%
\]

\end{corollary}

So the mirror toric Fano variety $Y^{\circ}$ as defined in Section
\ref{sec tropical mirror polarization} coincides with $\mathbb{P}\left(
\nabla\right)  $ as defined by Batyrev and Borisov.

\subsubsection{Construction of the mirror degeneration via tropical geometry}

We now describe the tropical subcomplex $T_{I_{0}}\left(  I\right)
\subset\partial\nabla_{I_{0}}\left(  I\right)  $. As before denote by
$\operatorname*{Strata}\nolimits_{\Delta}\left(  I_{0}\right)  $ the
subcomplex of toric strata of the boundary complex $\partial\Delta$ of
$\Delta$.

\begin{theorem}
\label{thm mirror map}For the degeneration defined in Lemma \ref{lem degen ci}
the map%
\[%
\begin{tabular}
[c]{ccccc}%
$\nabla_{I_{0}}\left(  I\right)  $ &  & $\nabla_{I_{0}}\left(  I\right)
^{\ast}$ &  & $\Delta$\\
$\cup$ &  & $\cup$ &  & $\cup$\\
$T_{I_{0}}\left(  I\right)  $ & $\rightarrow$ & $T_{I_{0}}\left(  I\right)
^{\ast}$ & $\rightarrow$ & $\operatorname*{Strata}\nolimits_{\Delta}\left(
I_{0}\right)  $\\
$F$ & $\mapsto$ & $F^{\ast}$ & $\mapsto$ & $\sum_{i=1}^{c}F^{\ast}\cap
\Delta_{i}$%
\end{tabular}
\
\]
is an inclusion reversing bijection.
\end{theorem}

Note, that this in particular shows that the complex $T_{I_{0}}\left(
I\right)  $ is dual to the sphere $\operatorname*{Strata}\nolimits_{\Delta
}\left(  I_{0}\right)  $.

\begin{proposition}
\label{prop TI0I}Let $\Sigma^{\circ}\left(  1\right)  =J_{1}^{\circ}%
\cup...\cup J_{c}^{\circ}$ be the nef partition corresponding to the
Batyrev-Borisov mirror and $I_{0}^{\circ}$ the associated monomial ideal as
defined in Section \ref{sec degen ci}. Then%
\[
\operatorname*{Strata}\nolimits_{\nabla}\left(  I_{0}^{\circ}\right)
=T_{I_{0}}\left(  I\right)
\]

\end{proposition}

From Theorem \ref{thm mirror map} and Proposition \ref{prop TI0I} we obtain
that the lattice points of the support of $\operatorname*{Strata}%
\nolimits_{\Delta}\left(  I_{0}\right)  ^{\ast}\subset\partial\Delta^{\ast}$
correspond to the first order deformations of $I_{0}^{\circ}$. As
$X_{0}^{\circ}$ is again a complete intersection, by the Koszul complex the
first order tropical mirror family is a global flat family:

\begin{theorem}
The tropical mirror degeneration of $\mathfrak{X}$ (as introduced in Section
\ref{sec first order mirror degeneration}) defines a flat family
$\mathfrak{X}^{\circ}\subset Y^{\circ}\times\operatorname*{Spec}%
\mathbb{C}\left[  s\right]  $ and this coincides with the degeneration
associated to the nef partition $\Sigma^{\circ}\left(  1\right)  =J_{1}%
^{\circ}\cup...\cup J_{c}^{\circ}$ (by Lemma \ref{lem degen ci}).
\end{theorem}

Of course, for a complete intersection it is sufficient to compute the mirror
special fiber $X_{0}^{\circ}$ (as this again is a complete intersection).
However the deformation data is necessary to describe the mirror degeneration
in the case of non-complete intersections like, e.g., Pfaffian varieties.

Note also, that we have reproduced the Batyrev-Borisov mirror without the
knowledge that it was a complete intersection, and without any non-trivial use
of convex geometry (i.e., aside from convex hulls), as the tropical mirror
construction directly obtains the relevant data $T_{I_{0}}\left(  I\right)
^{\ast}\subset\partial\nabla^{\ast}$ (see also the Algorithm formulated in
Section \ref{Sec rem}).

\section{Example and implementation\label{sec examples}}

We formulate the construction given in Section \ref{Sec Applications} in the
form of an algorithm as implemented by the author in the \textsc{Macaulay2}
\cite{GS} package \textsc{SRdeformations} \cite{SRdeformations} (though the
construction is not limited to complete intersections, but more complicated in
general). In the complete intersection case we just have to specify the
special fiber ideal $I_{0}$ and can obtain from that $\mathfrak{X}^{\circ}$
(and $\mathfrak{X}$).

\begin{algorithm}                      
\caption{Tropical mirror family}          
\label{alg trm}
\begin{algorithmic}[1]
\REQUIRE Monomial ideal $I_{0}$ corresponding to a nef complete intersection Calabi-Yau variety in a Gorenstein toric Fano variety $Y$.
\ENSURE The tropical mirror family $\mathfrak{X}^{\circ}
\subset
Y^{\circ}\times\operatorname*{Spec}\mathbb{C}[ s]$
\STATE Compute a torus invariant basis of $\operatorname*{Hom}\left(  I_{0}%
,S/I_{0}\right)  _{0}$.
\STATE Compute the convex hull $\nabla^{\ast}$ of the lattice monomials corresponding
to this basis.
\STATE Find the co-complex $T_{I_{0}}\left(  I\right)  ^{\ast}$ of tropical faces of
$\nabla^{\ast}$, i.e., those faces $F$ of $\nabla^{\ast}$ such that the ideal
\[
\phi_{F}\left(  I_{0}\right)  =\left\langle m_{0}+t\cdot%
{\textstyle\sum\nolimits_{\alpha\in F\cap M}}
c_{\alpha}\cdot\phi_{\alpha}\left(  m_{0}\right)  \mid m_{0}\in I_{0}%
\right\rangle \subset\mathbb{C}\left[  t\right]  \otimes S
\]
with generic coefficients $c_{\alpha}$ does not contain a monomial.
\STATE Via Equation \ref{equ I0mirror} obtain the ideal $I_{0}^{\circ}$ associated to $T_{I_{0}}\left(
I\right)  \subset\partial\nabla$.
\RETURN the first order tropical mirror family $\mathfrak{X}^{\circ}$ defined in Section
\ref{sec first order mirror degeneration}.
\end{algorithmic}
\end{algorithm}

\begin{example}
We treat the $K3$ surface given as the complete intersection of a quadric and
a cubic in $\mathbb{P}^{4}$ using the \textsc{Macaulay2} package:\smallskip

\noindent\texttt{i1: R = QQ[x\_0..x\_4];}

\noindent\texttt{i2: I0 = ideal(x\_0*x\_1,x\_2*x\_3*x\_4);}\smallskip

\noindent The Stanley-Reisner complex of $I_{0}$:\smallskip

\noindent\texttt{i3: C = idealToComplex I0;}

\noindent\texttt{o3: 2: x x x\hspace{0.3in}x x x\hspace{0.3in}x x
x\hspace{0.3in}x x x\hspace{0.3in}x x x\hspace{0.3in}x x x \vspace{-0.05in}}

\noindent\texttt{\hspace{0.3in}\texttt{\hspace{0.5in}}0 2 3\hspace{0.3in}1 2
3\hspace{0.3in}0 2 4\hspace{0.3in}1 2 4\hspace{0.3in}0 3 4\hspace{0.3in}1 3 4}

\noindent\texttt{\hspace{0.3in}complex of dim 2 embedded in dim 4 (printing
facets)}

\noindent\texttt{\hspace{0.3in}equidimensional, simplicial, F-vector
\{1,5,9,6,0,0\}}\smallskip

\noindent Computing $\nabla^{\ast}$ as the convex hull of the
deformations:\smallskip

\noindent\texttt{i4: NablaDual=PT1 C;}

\noindent\texttt{o4: 4: y y y y y y y y y y \vspace{-0.05in}}

\noindent\texttt{\hspace{0.3in}\texttt{\hspace{0.5in}}0 1 2 3 4 5 6 7 8 9}

\noindent\texttt{\hspace{0.3in}complex of dim 4 embedded in dim 4 (printing
facets)}

\noindent\texttt{\hspace{0.3in}equidimensional, non-simplicial, F-vector
\{1,10,24,25,11,1\}}\smallskip

\noindent Compute $T_{I_{0}}\left(  I\right)  ^{\ast}$ as a subcomplex of the
boundary of $\nabla^{\ast}$:\smallskip

\noindent\texttt{i5: TI0Dual = tropDef(C,NablaDual)}

\noindent\texttt{o5: 1: y y\hspace{0.3in}y y\hspace{0.3in}y y\hspace{0.3in}y
y\hspace{0.3in}y y \vspace{-0.05in}}

\noindent\texttt{\hspace{0.3in}\texttt{\hspace{0.5in}}0 4\hspace{0.3in}8
9\hspace{0.3in}3 7\hspace{0.3in}2 6\hspace{0.3in}1 5}

\noindent\texttt{\hspace{0.3in}co-complex of dim 1 embedded in dim 4 (printing
facets)}

\noindent\texttt{\hspace{0.3in}equidimensional, non-simplicial, F-vector
\{0,0,5,9,6,1\}}\smallskip

\noindent Dualize to obtain $T_{I_{0}}\left(  I\right)  $ as a subcomplex of
the boundary of $\nabla$:\smallskip

\noindent\texttt{i6: TI0 = dualize TI0Dual}

\noindent\texttt{o6: 2: v v v\hspace{0.3in}v v v v\hspace{0.3in}v v v
v\hspace{0.3in}v v v v\hspace{0.3in}v v v \vspace{-0.05in}}

\noindent\texttt{\hspace{0.3in}\texttt{\hspace{0.5in}}2 4 7\hspace{0.3in}2 4 8
9\hspace{0.3in}2 5 7 9\hspace{0.3in}4 5 7 8\hspace{0.3in}5 8 9}

\noindent\texttt{\hspace{0.3in}complex of dim 2 embedded in dim 4 (printing
facets)}

\noindent\texttt{\hspace{0.3in}equidimensional, non-simplicial, F-vector
\{1,6,9,5,0,0\}}\smallskip

\noindent The coordinates of the vertices of $\nabla_{I_{0}}\left(  I\right)
=\nabla$:\smallskip

\noindent\texttt{i7: transpose TI0.grading}

\noindent\texttt{o7:}%
\[
\left(
\begin{array}
[c]{ccccccccccc}%
1 & 0 & 1 & 0 & 1 & -1 & 0 & 1 & -1 & -1 & -1\\
0 & 1 & 1 & 0 & 0 & -1 & 0 & 0 & -1 & 0 & -1\\
0 & 0 & 0 & 1 & 1 & -1 & 0 & 0 & 0 & -1 & -1\\
0 & 0 & 0 & 0 & 0 & 0 & 1 & 1 & -1 & -1 & -1
\end{array}
\right)
\]

\noindent\texttt{i8: fvector C}

\noindent\texttt{o8: \{1, 5, 9, 6, 0, 0\}}

\noindent\texttt{i9: fvector B}

\noindent\texttt{o9: \{1, 6, 9, 5, 0, 0\}}

\noindent We observe that $T_{I_{0}}\left(  I\right)  $ and
$\operatorname*{Strata}\nolimits_{\Delta}\left(  I_{0}\right)  $ have mirror
dual $F$-vectors.
\end{example}

The code computing this example and others can be found in the documentation
of the package \textsc{SRdeformations} \cite{SRdeformations}.

\section{Remarks and further applications\label{Sec rem}}

The tropical mirror construction can also be applied for degenerations
$\mathfrak{X}$ of a non-generic complete intersection $X$ to $X_{0}$ (defined
by $I_{0}$) as long as $\nabla_{I_{0}}\left(  I\right)  $ will still be a
polytope, e.g., to handle subfamilies with prescribed singularities.

The construction is also applicable, e.g., to non-complete intersection
Gorenstein Calabi-Yau varieties of codimension $3$, indeed, handling
non-complete intersection cases is the main aim of the construction. In
particular, as will be treated in a separate paper, it can be used to
reproduce a known mirror construction for a Pfaffian Calabi-Yau by R\o dland
\cite{Rodland The Pfaffian CalabiYau its Mirror} and yields new examples of
non-complete intersection mirror pairs, see, e.g., \cite[Sec. 10]{B Mirror
symmetry and tropical geomety}.

\end{document}